\documentclass{article}

\usepackage{amssymb} 

\newtheorem{theorem}{Theorem}

\begin{document}
%
\title{The Analytical Solution of the Lag-Lead Compensator}

\author{Li Li\footnote{Corresponding author:
li-li@mail.tsinghua.edu.cn}, Zhengpeng Wu\\
\\
Department of Automation, Tsinghua University, Beijing, China
100084}

\maketitle

\begin{abstract}
In this paper, we first give the analytical solution of the general
lag-lead compensator design problem. Then, we show why a series of
more than $5$ phase-lead/phase-lag compensator cannot be solved
analytically using the Galois Theory.
\end{abstract}

\section{Introduction}
\label{sec:1}

The well known lag-lead compensator design problem is a typical
frequency controller design problem; see also the related
discussions in the textbooks listed in
\cite{DavisonChenPloenBernstein2007}. During the last four decades,
different design methods were proposed
\cite{Wakeland1974}-\cite{Wang2006}. The analytical design
procedures for single continuous phase-lag and phase-lead
compensator have been given in several literatures, e.g.
\cite{Wang2003}. An analytical solving procedure is constructed for
three-parameter lag-lead compensators in \cite{Wang2006}. But that
method cannot be directly applied to four-parameter cases. A
universal design chart based four-parameter lag-lead compensator
design method was proposed in \cite{YeungLee2000}. Though it makes
great progress to avoid manual graphical manipulations in design, it
is still a graph based approach and sometimes does not yield the
accurate solution. To our best knowledge, the analytical solution of
four-parameter or even more general lag-lead compensator remains
unsolved till now.

In this paper, we will first give the analytical solution of the
general lag-lead compensator. Then, we will show why a series of
more than $5$ phase-lead/phase-lag compensator usually cannot be
analytically determined using the Galois Theory.

\section{The Analytical Solution for the General Lag-Lead Compensator}
\label{sec:2}

In general, a $n$th-order lag-lead compensator ($n \ge 1$) can be
written as
\begin{equation}
\label{equ:1} G_c (s) = \frac{s^n + b_1 s^{n-1} + ... + b_n}{ s^n +
a_1 s^{n-1} + ... + a_n}
\end{equation}
\noindent where $a_i$, $b_i \in \mathbb{R}^+ \cup \{0\}$, for $i=1$,
..., $n$, due to the requirement of casual stability.

Substitute $s$ with $j \omega$, we get
\begin{equation}
\label{equ:2} \bar{G}_c (j \omega) = \frac{(j \omega)^n + b_1 (j
\omega)^{n-1} + ... + b_n}{(j \omega)^n + a_1 (j \omega)^{n-1} + ...
+ a_n}
\end{equation}

Usually, the dedicated performance requirements are given as several
pairs of gain and phase at certain frequencies. For the $k$th
performance requirement, we have
\begin{equation}
\label{equ:3} \bar{G}_c (j \omega_k) = \frac{(j \omega_k)^n + b_1 (j
\omega_k)^{n-1} + ... + b_n}{(j \omega_k)^n + a_1 (j \omega_k
)^{n-1} + ... + a_n} = g_k \cos(p_k) + g_k \sin(p_k) j
\end{equation}
\noindent where $g_k$ and $p_k$ are the corresponding gain and phase
at frequency $\omega_k$, for $k \in \mathbb{N}$.

Eq.(\ref{equ:3}) can be rewritten as
\begin{equation}
\label{equ:4} (j \omega_k)^n + b_1 (j \omega_k)^{n-1} + ... + b_n =
\left[ (j \omega_k)^n + a_1 (j \omega_k )^{n-1} + ... + a_n \right]
\left[ g_k \cos(p_k) + g_k \sin(p_k) j \right]
\end{equation}

I) If $n$ is an even integer satisfying $n = 2m$, $m \in
\mathbb{N}$. From Eq.(\ref{equ:4}), we can have
\begin{eqnarray}
\label{equ:5} & & (-1)^{m} \omega_k^{2m} + \sum_{q=1}^{m} (-1)^{m-q}
\omega_k^{2m-2q} b_{2q}
+ j \sum_{q=1}^{m} (-1)^{m-q} \omega_k^{2m-2q+1} b_{2q-1} \nonumber \\
& = & \left[ (-1)^{m} \omega_k^{2m} + \sum_{q=1}^{m} (-1)^{m-q}
\omega_k^{2m-2q} a_{2q} + j
\sum_{q=1}^{m} (-1)^{m-q} \omega_k^{2m-2q+1} a_{2q-1} \right] \left[ g_k \cos(p_k) + g_k \sin(p_k) j \right] \nonumber \\
& = & \left[ g_k \cos(p_k) \left((-1)^{m} \omega_k^{2m} +
\sum_{q=1}^{m} (-1)^{m-q} \omega_k^{2m-2q} a_{2q} \right) - g_k
\sin(p_k) \sum_{q=1}^{m} (-1)^{m-q}
\omega_k^{2m-2q+1} a_{2q-1} \right] \nonumber \\
& & + j \left[ g_k \sin(p_k) \left((-1)^{m} \omega_k^{2m} +
\sum_{q=1}^{m} (-1)^{m-q} \omega_k^{2m-2q} a_{2q} \right) + g_k
\cos(p_k) \sum_{q=1}^{m} (-1)^{m-q} \omega_k^{2m-2q+1} a_{2q-1}
\right] \nonumber \\
\end{eqnarray}
\noindent which finally leads to the following two linear equations
of $a_i$, $b_i$, for $i=1$, ..., $n$.

\begin{eqnarray}
\label{equ:6} & \sum_{q=1}^{m} (-1)^{m-q} \omega_k^{2m-2q} b_{2q} -
g_k \cos(p_k) \sum_{q=1}^{m} (-1)^{m-q} \omega_k^{2m-2q} a_{2q} \nonumber \\
& + g_k \sin(p_k) \sum_{q=1}^{m} (-1)^{m-q} \omega_k^{2m-2q+1}
a_{2q-1} = - (-1)^{m} \omega_k^{2m} + g_k \cos(p_k) (-1)^{m}
\omega_k^{2m}
\end{eqnarray}
%
\begin{eqnarray}
\label{equ:7} & \sum_{q=1}^{m-1} (-1)^{m-q} \omega_k^{2m-2q+1}
b_{2q-1} - g_k \sin(p_k) \sum_{q=1}^{m} (-1)^{m-q} \omega_k^{2m-2q}
a_{2q} \nonumber \\
& - g_k \cos(p_k) \sum_{q=1}^{m} (-1)^{m-q} \omega_k^{2m-2q+1}
a_{2q-1} = g_k \sin(p_k) (-1)^{m} \omega_k^{2m}
\end{eqnarray}

II) Similarly, if $n$ is an odd integer satisfying $n = 2m-1$, $m
\in \mathbb{N}$, the $k$th performance requirement will also lead to
two linear equations of $a_i$, $b_i$, for $i=1$, ..., $n$.

In the rest of this paper, we will call $r$ performance requirement
pairs ($g_k$, $p_k$, $\omega_k$), $k=1$, ..., $r$ are feasible, if
they lead to a $2r$ consistent and linearly independent
(irreducible) equation set defined as (\ref{equ:6})-(\ref{equ:7}).
As a result, we can reach the following conclusion.

\indent

\begin{theorem} Suppose we have $r$ feasible performance
requirement pairs ($g_k$, $p_k$, $\omega_k$), $k=1$, ..., $r$. If $r
< n$, we may have infinite possible solutions of this compensator.
If $r > n$, we cannot find a feasible solution of this compensator.
If $r = n$, we can formulate a $2n$ consistent and linearly
independent linear equation set for these $2n$ unknown parameters
$a_i$, $b_i$, for $i=1$, ..., $n$. Thus, we can get the analytical
solution of this lag-lead compensator directly by solving this
linear equations set (e.g. using Cramer's rule).
\end{theorem}

\indent

It is easy to prove that the analytical solving methods of
phase-lag/phase-lead and three-parameter lag-lead compensator design
problem proposed in \cite{Wang2003}-\cite{Wang2006} are indeed
special cases of the above method.

\section{Further Discussions}
\label{sec:3}

There are two interesting questions concerning the lag-lead
compensator design problems. The first question is

\indent

\textbf{Question 1}: Determine whether a set of performance
requirement pairs ($g_k$, $p_k$, $\omega_k$), $k=1$, ..., $n$ is
feasible for a $n$th-order lag-lead compensator.

\indent

From the above discussion, we can see that a set of $n$ performance
requirement pairs is feasible unless they lead to $2n$ consistent
and linearly independent. Moreover, it is often required the
lag-lead compensator to be casual stable. Thus, we need to check the
algebraic stability criterion for the following equation
\begin{equation}
\label{equ:8} s^n + a_1 s^{n-1} + ... + a_n = 0
\end{equation}
\noindent after obtaining $a_i$, $i=1$, ..., $n$.

The necessary and sufficient algebra stability criterion for
Eq.($\ref{equ:8}$) is hard to find. However, we can apply
Routh-Hurwitz stability criterion which is necessary and frequently
sufficient. Since readers are familiar with this issue, we will not
discuss the details.

The second question is

\indent

\textbf{Question 2}: Determine whether we find a series of $n$
phase-lead/phase-lag compensator connected as
\begin{equation}
\label{equ:9} G_c (s) = \frac{s + d_1}{s + c_1} \cdot \frac{s +
d_2}{s + c_2} \cdot ... \cdot \frac{s + d_n}{s + c_n}
\end{equation}
\noindent which can satisfy a set of performance requirement pairs
($g_k$, $p_k$, $\omega_k$), $k=1$, ..., $n$. Here, $c_i$, $d_i \in
\mathbb{R}$, for $i=1$, ..., $n$.

\indent

From the above discussion, if this set of performance requirement
pairs ($g_k$, $p_k$, $\omega_k$), $k=1$, ..., $n$, is feasible, we
have
\begin{equation}
\label{equ:10} G_c (s) = \frac{(s + d_1)(s + d_2) ... (s + d_n)}{(s
+ c_1)(s + c_2) ... (s + c_n)} = \frac{s^n + b_1 s^{n-1} + ... +
b_n}{s^n + a_1 s^{n-1} + ... + a_n}
\end{equation}
\noindent where $a_i$, $b_i$, $i=1$, ..., $n$ are calculated from
the selected performance requirements using the above method.

Thus, \textit{Question 2} is equal to finding the roots of $s^n +
a_1 s^{n-1} + ... + a_n = 0$ and $s^n + b_1 s^{n-1} + ... + b_n =
0$.

Based on the well known Galois theory \cite{Lang2002}-
\cite{Stewart2003}, we can always find the the analytical solution
of $c_i$, $d_i$, $i=1$, ..., $n$ for $n \in \{1, 2, 3, 4 \}$. But
generally, we cannot find the analytical solution for $n \ge 5$.



\begin{thebibliography}{1}


\bibitem{DavisonChenPloenBernstein2007}
D. E. Davison, J. Chen, O. R. Ploen, and D. S. Bernstein, ``What is
your favorite book on classical control?'' \emph{IEEE Control
Systems Magazine}, vol. 27, no. 3, pp. 89-99, 2007.

\bibitem{Wakeland1974}
W. R. Wakeland, ``Bode compensator design,'' \emph{IEEE Transactions
on Automatic Control}, vol. 21, no. 5, pp. 771-773, 1974.

\bibitem{Mitchell1977}
J. R. Mitchell, ``Comments on `Bode compensator design',''
\emph{IEEE Transactions on Automatic Control}, vol. 22, no. 5, pp.
869-870, 1977.

\bibitem{YeungChaidDinh1995}
K. S. Yeung, K. Q. Chaid and T. X. Dinh, ``Bode design charts for
continuous-time and discrete-time compensators,'' \emph{IEEE
Transactions on Education}, vol. 38, no. 2, pp. 252-257, 1995.

\bibitem{YeungWongChen1998}
K. S. Yeung, K. W. Wong and K.-L. Chen, ``A non-trial-and-error
method for lag-lead compensator design,'' \emph{IEEE Transactions on
Education}, vol. 41, no. 1, pp. 74-80, 1998.

\bibitem{YeungLee2000}
K. S. Yeung and K. H. Lee, ``A universal design chart for linear
time-invariant continuous-time and discrete-time compensators,''
\emph{IEEE Transactions on Education}, vol. 43, no. 3, pp. 309-315,
2000.

\bibitem{Calleja2003}
H. Calleja, ``An approach to amplifier frequency compensation,''
\emph{IEEE Transactions on Education}, vol. 46, no. 1, pp. 43-49,
2003.

\bibitem{Wang2003}
F.-Y. Wang, ``The exact and unique solution for phase-lead and
phase-lag compensation,'' \emph{IEEE Transactions on Education},
vol. 46, no. 2, pp. 258-262, 2003.

\bibitem{Wang2006}
F.-Y. Wang, ``A new non-trial-and-error method for lag-lead
compensator design: A special case,'' \emph{International Journal of
Intelligent Control and Systems}, vol. 11, no. 1, pp. 69-76, 2006.


\bibitem{Lang2002}
S. Lang, \emph{Algebra}, Graduate Texts in Mathematics, vol. 211,
3rd edition, Springer-Verlag, New York, NY, 2002.

\bibitem{Edwards1984}
H. M. Edwards, \emph{Galois Theory}, Graduate Texts in Mathematics ,
vol. 101, Springer, New York, NY, 1984.

\bibitem{Stewart2003}
I. Stewart, \emph{Galois Theory}, 3rd edition, Chapman and Hall/CRC,
Boca Raton, FL, 2003.

\end{thebibliography}
\end{document}